\documentclass[reqno,12pt]{amsart}
\usepackage{amsmath,epsfig}

\newtheorem{theorem}{Theorem}
\newtheorem{lemma}{Lemma}
\newtheorem{proposition}{Proposition}
\newtheorem{remark}{Remark}
\newcommand{\Cov}{\operatorname{Cov}}
\newcommand{\Var}{\operatorname{Var}}
\newcommand{\Per}{\operatorname{Per}}
\newcommand{\VOR}{\operatorname{Vor}}
\newcommand{\Vol}{\operatorname{Vol}}
\newcommand{\la}{{\lambda}}
\newcommand{\XX}{{\mathcal X}}
\newcommand{\X}{{\mathcal X}}

\def\R{\mathbb{R}}

\def\0{{\bf 0}}

\begin{document}
\setlength{\baselineskip}{1.3\baselineskip}

\title{\bf Process level moderate deviations for stabilizing functionals}

\author{P. Eichelsbacher and T. Schreiber}

\maketitle

\footnotetext{{\em American Mathematical Society 2000 subject
classifications.} Primary 60F05, Secondary 60D05 } 

\begin{abstract}
Functionals of spatial point process often satisfy a weak spatial dependence
condition known as {\it stabilization}. In this paper we prove process level 
moderate deviation principles (MDP) for such functionals, which is
a level-3 result for empirical point fields as well as a level-2 result
for empirical point measures. The level-3 rate function coincides with
the so-called specific information. We show that the general result
can be applied to prove MDPs for various particular functionals,
including random sequential packing, birth-growth models, germ-grain
models and nearest neighbor graphs.  
\end{abstract}

\section{Introduction and main results}
 \subsection{Terminology}
 Consider a real-valued translation invariant functional
 $\xi(x,\sigma)$ defined on all pairs $(x,\sigma),$
 where $x \in {\Bbb R}^d$ and $\sigma$ is a finite point
 configuration in ${\Bbb R}^d$ containing $x.$ Moreover,
 for $x \not \in \sigma$ write $\xi(x,\sigma) :=
 \xi(x,\sigma \cup \{ x \}).$ Let ${\mathcal P}$ be a homogeneous
 Poisson point process on ${\Bbb R}^d,$ with a certain intensity
 $\tau > 0$ to remain fixed throughout the paper, and denote by
 $\Pi$ the distribution of ${\mathcal P}$ on the space $\Sigma$ of
 locally finite point configurations in ${\Bbb R}^d.$ For
 formal completeness we represent the space $\Sigma$ as
 the set of all locally finite and simple (all atoms 
 of mass 1) counting measures $\sigma$ on ${\Bbb R}^d,$
 endowed with the $\sigma$-field ${\mathcal F}$ generated by the 
 mappings $\Sigma \ni \sigma \mapsto \sigma(A)$ for
 all bounded Borel $A \subseteq {\Bbb R}^d.$ 

 One crucial assumption imposed on $\xi$ throughout this paper is the so-called 
 exponential stabilization, see \cite{BY,Pe2,Pe3,PY1}.
 We say that $\xi$ is stabilizing (at intensity $\tau$)
 if for each $x \in {\Bbb R}^d$ there exists an a.s. finite random variable
 $R(x) := R^{\xi}(x,{\mathcal P})$ (a radius of stabilization) and
 $\xi_{\infty}(x) := \xi_{\infty}(x,{\mathcal P})$ (the limit of $\xi$)
 such that, with probability one, $\xi(x,({\mathcal P}
 \cap B_{R(x)}(x)) \cup \sigma) = \xi_{\infty}(x)$
 for all locally finite $\sigma \subseteq {\Bbb R}^d
 \setminus B_{R(x)}(x).$ More generally, for a locally
 finite point configuration $\sigma \subseteq {\Bbb R}^d$ we
 consider the stabilization radius $R(x) := R^{\xi}(x,\sigma)$
 of $\xi$ at $x$ in $\sigma$ defined so that
 $\xi(x, (\sigma \cap B_{R(x)}(x)) \cup \sigma')$ takes
 the same value for all locally finite $\sigma' \subseteq
 {\Bbb R}^d \setminus B_{R(x)}(x).$ We put $R^{\xi}(x,\sigma) 
 := +\infty$ if this does not hold for any finite $R(x).$ 
 Our exponential stabilization requirement means that 
 $\xi$ is stabilising (at the intensity $\tau$ fixed
 throughout the paper) and at each point the stabilisation
 radius exhibits exponentially decaying tail, i.e.
 
   {\bf (E)} There exists $c > 0$ such that,
        for $r$ large enough, 
        $$ {\Bbb P}(R(x) > r) \leq \exp(-c r). $$

In \cite{BY,Pe2,Pe3,PY1} a lot of examples of stabilizing functionals are discussed.
In Section 2 we will focus on random sequential packing models, birth-growth models, 
germ-grain models and nearest neighbor graphs.

 Under the stabilization condition as stated above, the Poisson
 point process ${\mathcal P}$ with probability one takes its values
 in the space $\Sigma^{\xi} \subseteq \Sigma,$ defined to
 consist of all configurations for which the value of $\xi$
 can be uniquely determined at each configuration point. 
 Thus, in order to avoid unnecessary formal subtleties,
 we simply extend the functional $\xi$ in some artificial
 way, say by setting $\xi(x,\sigma) := 0$ if $R^{\xi}(x,\sigma)
 = +\infty.$ Since this happens with probability $0$ if
 $\sigma$ is given by ${\mathcal P},$ this extension does
 not affect our results while guaranteeing that $\Sigma^{\xi}
 = \Sigma.$ For a configuration $\sigma \in \Sigma$ let 
 $\xi[\sigma]$ be its $\xi$-marked version, where
 each point $x \in \sigma$ is marked with the
 corresponding value $\xi(x,\sigma).$ In particular,
 $\xi[\Sigma]$ is the space of all possible
 $\xi$-marked point configurations. We formally 
 represent $\xi[\Sigma]$ as the space of simple point measures of
 ${\Bbb R}^d \times {\Bbb R}$ and we endow it with the $\sigma$-field
 $\hat{\mathcal F}$ generated by the mappings $\xi[\Sigma] \ni \hat{\sigma} 
 \mapsto \hat{\sigma}(A_1 \times A_2)$ for all bounded Borel $A_1 \subseteq
 {\Bbb R}^d,\; A_2 \subseteq {\Bbb R}.$  

 For each Borel measurable region $A \subseteq {\Bbb R}^d$ consider the
 $\sigma$-field ${\mathcal F}_A \subseteq {\mathcal F}$ generated by the
 mappings $\Sigma \ni \sigma \mapsto \sigma(B),\; B \subseteq
 A,\; B$ bounded and measurable. Define also $\hat{\mathcal F}_A$
 to be the $\sigma$-field generated by the mappings $\xi[\Sigma] \ni
 \hat\sigma \mapsto \hat\sigma(B_1 \times B_2)$ with $B_1$ ranging over 
 bounded Borel subsets of $A$ and with bounded Borel measurable $B_2
 \subseteq {\Bbb R}.$ We shall write $\Pi_A$ for the restriction of $\Pi$
 to ${\mathcal F}_A.$ We say that a function
 $\Phi : \Sigma \to {\Bbb R}$ is local if it is measurable with
 respect to ${\mathcal F}_A$ for some bounded $A.$ Likewise,
 $\hat{\Phi} : \xi[\Sigma] \to {\Bbb R}$ is local if it is 
 measurable with respect to $\hat{\mathcal F}_A$ for some
 bounded $A.$ Consider the space $B_{loc}(\Sigma)$ consisting
 of all the bounded local functions on $\Sigma$ with the
 topology determined by the convergence: $\Phi_n \to \Phi$
 as $n \to \infty$ iff $||\Phi_n-\Phi||_{\infty} :=
 \sup_{\sigma \in \Sigma} |\Phi_n(\sigma) - \Phi(\sigma)|
 \rightarrow_{n\to\infty} 0$ and there exists bounded 
 $A \subseteq {\Bbb R}^d$ such that $\Phi,\Phi_1,\Phi_2,
 \ldots$ are all ${\mathcal F}_A$-measurable. The definition
 of $B_{loc}(\xi[\Sigma])$ is completely analogous. 
 
 We say that a set function $\Theta : {\mathcal F} \to {\Bbb R}$ is a
 signed local measure on $\Sigma$ iff $\Theta(\bigcup_{i=1}^{\infty} S_i)
 = \sum_{i=1}^{\infty} \Theta(S_i)$ with the RHS series absolutely convergent,
 whenever $S_i$ are pairwise disjoint and all $S_i$ are ${\mathcal F}_A$-measurable
 for some bounded $A \subseteq {\Bbb R}^d.$ Denote by
 ${\mathcal M}^{0,\theta}_{loc}(\Sigma)$
 the space of all translation invariant signed local measures
 on $\Sigma$ with total mass $\int_{\Sigma} 1 d \Theta = 0$ (and
 hence referred to as null-measures in the sequel) endowed with the topology
 ${\mathcal T}$ taken to be the weakest one which makes continuous
 the mappings $\Theta \mapsto \langle \Phi, \Theta \rangle
 := \int_{\Sigma} \Phi d\Theta$ 
 for all $\Phi \in B_{loc}(\Sigma).$ Observe that the mapping
 $\Phi \mapsto \langle \Phi, \Theta \rangle$ is continuous
 in $B_{loc}(\Sigma)$ for each $\Theta \in
 {\mathcal M}^{0,\theta}_{loc}(\Sigma).$ Moreover, it is
 clear that $B_{loc}(\Sigma)$ is the topological dual of
 $({\mathcal M}^{0,\theta}_{loc}(\Sigma),{\mathcal T}).$
 In a completely
 analogous way we define the space ${\mathcal M}^{0,\theta}_{loc}
 (\xi[\Sigma])$ of translation invariant signed local null-measures
 on $\xi[\Sigma],$ endowed with the corresponding
 $B_{loc}(\xi[\Sigma])$-weak topology $\hat{\mathcal T}^{\xi}.$
 Consider now the natural mark-forgetting mapping  
 $\pi : \xi[\Sigma] \to \Sigma$ and observe that 
 we have $\hat{\Theta} \circ \pi^{-1} \in
 {\mathcal M}^{0,\theta}_{loc}(\Sigma)$ for $\hat{\Theta}
 \in {\mathcal M}^{0,\theta}_{loc}(\xi[\Sigma]).$ 
 Define 
 $$ {\mathcal M}^{0,\theta}_{loc,\xi}(\Sigma) :=
   \left\{ \Theta \in {\mathcal M}^{0,\theta}_{loc}(\Sigma)\; | \;
    \exists \hat{\Theta} \in {\mathcal M}^{0,\theta}_{loc}(\xi[\Sigma]) :
    \Theta = \hat{\Theta} \circ \pi^{-1} \right\} $$
 and endow ${\mathcal M}^{0,\theta}_{loc,\xi}(\Sigma)$ with the
 topology ${\mathcal T}^{\xi} := \pi(\hat{\mathcal T}^{\xi}).$ Note that
 ${\mathcal M}^{0,\theta}_{loc,\xi}(\Sigma) \subseteq {\mathcal M}^{0,\theta}_{loc}(\Sigma)$
 where the equality may but does not have to occur. Indeed, for $\Theta \in 
 {\mathcal M}^{0,\theta}_{loc,\xi}(\Sigma)$ the existence of $\hat{\Theta} \in
 {\mathcal M}^{0,\theta}_{loc}(\xi[\Sigma])$ with $\Theta = \hat\Theta \circ \pi^{-1}$
 may imply a version of $\sigma$-additivity stronger than just local whenever 
 $\xi$ itself is not a local functional. For similar reasons, the topology 
 ${\mathcal T}^{\xi}$ is stronger or equal to the topology induced by the
 inclusion ${\mathcal M}^{0,\theta}_{loc,\xi}(\Sigma) \subseteq
 {\mathcal M}^{0,\theta}_{loc}(\Sigma).$ Observe also that ${\mathcal T}^{\xi}$
 can be equivalently characterized as the weakest topology to make continuous
 the mappings
 ${\mathcal M}^{0,\theta}_{loc,\xi}(\Sigma) \ni \Theta \mapsto
  [\hat{\Phi},\Theta]_{\xi} := \langle \hat{\Phi}, \xi[\Theta] \rangle$
 for all $\hat{\Phi} \in B_{loc}(\xi[\Sigma]).$ Clearly,
 $\hat{\Phi} \mapsto [\hat{\Phi},\Theta]_{\xi}$ is $B_{loc}(\xi[\Sigma])$-continuous
 for each $\Theta$ in ${\mathcal M}^{0,\theta}_{loc,\xi}(\Sigma).$ Moreover, it
 is easily seen that $B_{loc}(\xi[\Sigma])$ can be regarded as the topological dual
 for $({\mathcal M}^{0,\theta}_{loc,\xi}(\Sigma),{\mathcal T}^{\xi})$
 with respect to the duality $[\cdot,\cdot]_{\xi}.$

 For $\lambda > 0,$ let $Q_{\lambda}$ be the cube
 of volume $\lambda$ centered at $0,$  i.e. $Q_{\lambda} =
 [-\sqrt[d]{\lambda} \slash 2, \sqrt[d]{\lambda} \slash 2]^d.$
 For a finite point configuration $\sigma \subseteq Q_{\lambda}$
 we define the {\it empirical point field} 
 \begin{equation}\label{OKRESLENIER}
  \psi^{\xi}_{\lambda}(\sigma) := \frac{1}{\lambda}
  \int_{Q_{\lambda}} \delta_{\xi[\tau_x \Per_{\lambda}(\sigma)]} dx,
 \end{equation}
 where $\tau_x y := y-x$ stands for the usual shift operator,
 while $\Per_{\lambda}(\sigma)$ is the configuration arising by
 periodically copying $\sigma$ on disjoint translates of $Q_{\lambda},$
 i.e.
 $\Per_{\lambda}(\sigma) := \bigcup_{i \in {\Bbb Z}^d}
  \tau_{\sqrt[d]{\lambda}i} \sigma.$
 In other words, the empirical process arises as a probability measure
 on the space $\xi[\Sigma]$ of marked point configurations, by normalized
 integration over $x \in Q_{\lambda}$ of unit masses concentrated at
 $\xi$-marked $\sqrt[d]{\lambda}$-periodized versions of $\sigma$
 shifted by $x.$ It is clear that $\psi^{\xi}_{\lambda}(\sigma)$ is
 a translation invariant measure. Throughout this paper we focus
 our interest on empirical point fields generated by the Poisson
 point process ${\mathcal P}$
 \begin{equation}\label{EMPIPSI}
  \Psi^{\xi}_{\lambda} := 
  \psi^{\xi}_{\lambda}({\mathcal P} \cap Q_{\lambda}).
 \end{equation}
 We consider also the centered versions
 $$ \bar{\Psi}^{\xi}_{\lambda} := 
    \Psi^{\xi}_{\lambda} - {\Bbb E} \Psi^{\xi}_{\lambda}. $$
 Observe that $\bar{\Psi}^{\xi}_{\lambda}$ is a 
 ${\mathcal M}^{0,\theta}_{loc}(\xi[\Sigma])$-valued random
 element and that we almost surely have
 $\pi(\bar{\Psi}^{\xi}_{\lambda}) \in
  {\mathcal M}^{0,\theta}_{loc,\xi}(\Sigma).$ 
 It can be shown that the following law of large numbers holds almost surely in
 ${\mathcal T}^{\xi}$ topology 
 $$ \lim_{\lambda \to \infty} \Psi^{\xi}_{\lambda} 
   = \lim_{\lambda\to\infty}  {\Bbb E} \Psi^{\xi}_{\lambda} =
   \xi[\Pi], $$
 this is a consequence of the exponential stabilization and
 we refer the reader to \cite{PY4} for details. 
 The main goal of this paper is to establish a process level
 (level-3) and empirical measure level (level-2) moderate deviation
 principle (MDP) for $\bar{\Psi}^{\xi}_{\lambda}$ under the assumption
 that the stabilizing $\xi$ satisfies a level-1 moderate deviation
 principle, as made precise below. 
 The rate function of this MDP turns out to admit representation
 in terms of the {\it specific information} functional
 $I(\cdot|\Pi)$ defined for a local null-measure $\Theta \in
 {\mathcal M}^{0,\theta}_{loc}(\Sigma)$ by  
 \begin{equation}\label{GESTINF}
   I(\Theta | \Pi) := \frac{1}{2} \lim_{\lambda \to \infty} 
     \frac{1}{\lambda} \int_{\Sigma_{Q_{\lambda}}}
     \left( \frac{d\Theta_{|Q_{\lambda}}}{d \Pi_{Q_{\lambda}}}\right)^2
      d\Pi_{Q_{\lambda}}
  \end{equation}
  if $\Theta \ll_{loc} \Pi$ and $I(\Theta) := + \infty$ otherwise.
  Note that the local absolute continuity requirement stated
  with $\ll_{loc}$ above means simply that $\Theta_{|Q_{\lambda}}
  \ll \Pi_{Q_{\lambda}}$ for all $\lambda,$
  with the ${|Q_{\lambda}}$ operation standing
  for the restriction of its argument measure to ${\mathcal F}_{Q_{\lambda}}.$  
  The existence of the limit in (\ref{GESTINF}) will be established
  in Lemma \ref{ISTNIENIE}, further properties of the specific
  information functional will be discussed in Section \ref{GESTINFWL}.  

\subsection{Process level moderate deviation principles}
 To proceed with the statement of the moderate deviation principle
 for $\bar{\Psi}^{\xi}_{\lambda},$ we let $\alpha_{\lambda}$ be such
 that $\alpha_{\lambda} \to \infty$ and $\alpha_{\lambda}
 \lambda^{-1 \slash 2} \to 0.$ We say that a family of probability
 measures $(\mu_{\varepsilon})_{\varepsilon >0}$, on some
 topological space $\mathcal Y$ obeys a large deviation principle (LDP)
 with speed $\varepsilon$ and good rate function $I(\cdot): {\mathcal
 Y} \to \R^+_0 \cup \{ +\infty \}$ if
 \begin{itemize}
 \item $I$ is lower semi-continuous and has compact level sets
 $N_L:=\{x\in {\mathcal Y}: I(x) \le L\}$, for every $L \in  [0,
 \infty)$.
 \item
 For every open set $G \subseteq \mathcal Y$ it holds
 \begin{equation}
 \liminf_{\varepsilon \to 0} \varepsilon \log \mu_{\varepsilon}(G)\ge -\inf_{x\in G} I(x).
 \end{equation}
 \item
 For every closed set $ A\subseteq \mathcal Y$ it holds
 \begin{equation}
 \limsup_{\varepsilon \to 0} \varepsilon \log \mu_{\varepsilon}(A)\le -\inf_{x\in A} I(x).
 \end{equation}
 \end{itemize}
 Similarly we will say that a family of random variables
 $(Y_{\varepsilon})_{\varepsilon >0}$ with topological state space
 $\mathcal Y$ obeys a large deviation principle with speed
 $\varepsilon$ and good rate function $I(\cdot): {\mathcal Y} \to
 \R^+_0 \cup \{ +\infty \}$ if the sequence of their distributions
 does. Formally a moderate deviation principle is nothing but an
 LDP. However, we will speak about a moderate deviation principle
 (MDP) for a sequence of random variables, whenever the scaling of
 the corresponding random variables is between that of an ordinary
 law of large numbers and that of a central limit theorem.
 \medskip

 Below, we shall assume that $\xi$ is a bounded exponentially stabilizing
 functional, as required in {\bf (E)}. From the results and methods of
 Section 4.3 in \cite{BY} it follows that
 \begin{proposition}\label{VarianceLemma}
        For each $\hat\Phi \in B_{loc}(\xi[\Sigma])$ there exists the limit 
        \begin{equation}\label{VAAAR}
           V[\xi;\hat{\Phi}] := \lim_{\lambda\to\infty} \lambda
           \Var\left(\langle \hat{\Phi}, \bar \Psi^{\xi}_{\lambda} \rangle \right)
        \end{equation}
        providing the infinite-volume {\it variance density} for
        $\langle \hat{\Phi}, \bar \Psi^{\xi}_{\lambda} \rangle.$
        Moreover, for each $R>0$ we have  
        \begin{equation}\label{VAAAR2}
         V_R[\xi] := \sup_{||\hat\Phi||_{\infty} \leq 1, \; D(\hat\Phi) \leq R}
         V[\xi;\hat\Phi] < +\infty,
        \end{equation}
        where $D(\hat{\Phi})$ stands for the infimum of $r>0$ such that
        $\hat{\Phi}$ is ${\mathcal F}_{B(0,r)}$-measurable.    
       \end{proposition}
  Note that we multiply rather than dividing by $\lambda$ 
  in (\ref{VAAAR}) because of the normalization for $\Psi^{\xi}_{\lambda}$
  being already present in (\ref{OKRESLENIER}) and (\ref{EMPIPSI}).  
  Further, we impose on $\xi$ the following additional condition
 
   
  {\bf (L)} For the log-Laplace functional  
      \begin{equation}\label{PRZEKSZTLAPL}
       \Lambda^{\xi}_{\lambda;\alpha_{\lambda}}(\hat{\Phi}) :=
       \frac{1}{\alpha_{\lambda}^2} \log {\Bbb E} \exp
       \left( \alpha_{\lambda} \lambda^{1 \slash 2}
       \langle \hat{\Phi}, \bar{\Psi}^{\xi}_{\lambda} \rangle \right),\;\;
       \hat{\Phi} \in B_{loc}(\xi[\Sigma]) 
     \end{equation}

     we have
     \begin{equation}\label{LAPLL}
       \lim_{\lambda\to\infty} \Lambda^{\xi}_{\lambda;\alpha_{\lambda}}(\hat{\Phi})
       = \frac{1}{2} V[\xi;\hat{\Phi}].
     \end{equation}

  In fact, this condition is a usual ingredient needed to establish the moderate
  deviation principle for $\langle \hat{\Phi}, \bar \Psi^{\xi}_{\lambda} \rangle$
  with rate function ${\mathbb R} \ni t \mapsto t^2 \slash (2 V[\xi;\hat\Phi])$
  by an application of the G\"artner-Ellis theorem [Theorem 2.3.6 in \cite{DZ}],
  see \cite{BESY}. In a number of cases the exponential stabilization seems
  to be enough to guarantee {\bf (L)}, see \cite{BESY} and Section \ref{EXAMPLES}
  below, however at present we do not know if the boundedness and exponential
  stabilization do imply the condition {\bf (L)} in general.  

  The following process-level moderate deviation theorem is the first main result
  of our paper. 
 
\begin{theorem}\label{GLOWNE}
 For a bounded functional $\xi$ for which 
 the conditions {\bf (E)} and ${\bf (L)}$ hold, the family
 $(\alpha_{\lambda}^{-1} \lambda^{1 \slash 2}
 \bar{\Psi}^{\xi}_{\lambda})_{\lambda}$ satisfies the moderate
 deviation principle on the space $({\mathcal M}^{0,\theta}_{loc}(\xi[\Sigma]),
 \hat{\mathcal T}^{\xi})$ with speed $\alpha_{\lambda}^2$ and with the
 good rate function $I^{\xi}(\cdot|\xi[\Pi])$ defined for
 $\hat{\Theta} \in {\mathcal M}^{0,\theta}_{loc}(\xi[\Sigma])$
 to be $I(\Theta|\Pi)$ if $\hat{\Theta} = \xi[\Theta]$ and
 $+\infty$ otherwise.
\end{theorem}
  
 At this point, it is very natural to compare our Theorem \ref{GLOWNE}
 for stabilizing functionals with the corresponding process level
 large deviation principles for Gibbs measures, see F\"ollmer
 \& Orey \cite{FO}, Olla \cite{OL} and Georgii \cite{GEO}, where
 the rate function was given in terms of the relative entropy
 density. In fact, the specific information can be roughly
 interpreted as the (halved) second derivative of the relative
 entropy density $h(\cdot|\cdot)$ at the equilibrium measure $\Pi$
 in that, vaguely, $h(\Pi+\delta \Theta|\Pi) \approx \delta^2
 I(\Theta|\Pi) + o(\delta^2).$ Of course in such formulation
 this imprecise formula can be given a definite meaning only
 at the level of finite volume approximations of $h(\cdot|\Pi)$
 and $I(\cdot|\Pi),$ yet it provides an intuition that our 
 MDP could be regarded as a {\it local version} of the process
 level LDP, {\it differentiated} at equilibrium. To the best of our knowledge 
 there is no moderate deviations result on process level in the literature.


\subsection{Empirical measure level moderate deviation principles}

 Usually as a consequence from the process level (level-3) MDP one obtains its empirical
 measure level (level-2) counterpart, which is proven via the contraction principle.
 In our present context we prefer, however, to establish the level-2 MDP directly,
 thus avoiding certain topological intricacies and getting a simpler formula
 for the rate function, still in a variational rather than explicit form though.

 Let us consider the {\it empirical point process}
 \begin{equation}
  Z_{\lambda}^{\xi} := \frac{1}{\lambda}\sum_{x \in \mathcal P_{\tau} \cap Q_{\lambda}} 
  \delta_{\xi(x,\Per_{\lambda}({\mathcal P} \cap Q_{\lambda}))}.
 \end{equation}
 and its centered version $\bar{Z}_{\lambda}^{\xi} := Z_{\lambda}^{\xi} - 
 {\Bbb E} Z_{\lambda}^{\xi}$.
 Moreover let us denote by ${\mathcal M}(\R)$ the real vector space of finite
 variation signed measures on $\R$. Equip ${\mathcal M}(\R)$ with
 the {\it weak topology} generated by the sets $\{U_{f,x,\delta}, f
 \in {\mathcal C}_b(\R), x \in \R, \delta > 0 \}$, where ${\mathcal C}_b(\R)$ 
 is the space of bounded continuous functions on $\R$ and with 
 $$
  U_{f,x,\delta} := \{ \nu \in {\mathcal M}(\R),\; | \langle f, \nu \rangle - x | < \delta \}.
 $$
 The Borel-$\sigma$-field generated by the weak topology is denoted by ${\mathcal B}$.
 It is well known, that since the collection of linear functionals
 $\{\nu \mapsto \langle f, \nu \rangle: f \in {\mathcal C}_b(\R) \}$
 is separating in ${\mathcal M}(\R)$, this topology makes ${\mathcal M}(\R)$ into
 a locally convex, Hausdorff topological vector space, whose topological dual is the preceding
 collection, hereafter identified with ${\mathcal C}_b(\R)$. 
 \medskip

 In analogy with the corresponding results for process level objects, we require 
 that $\xi$ satisfy the exponential stabilization condition {\bf (E)}.
 Under this conditions, using the results of \cite{PY4}, we get the following
 almost sure law of large numbers in the ${\mathcal C}_b(\R)$-weak topology
 \begin{equation}\label{NUZBIE}
  \lim_{\lambda \to \infty} Z^{\xi}_{\lambda} = \lim_{\lambda \to \infty} 
   {\Bbb E} Z^{\xi}_{\lambda} = \tau \nu[\xi],
 \end{equation}
 where, for Borel $B \subseteq \R,$ 
 \begin{equation}\label{DEFNU}
  \nu[\xi](B) := {\Bbb P}( \xi({\bf 0},{\mathcal P}) \in B),
 \end{equation}
 that is to say, $\nu[\xi]$ is the law of $\xi({\bf 0},{\mathcal P})$ on $\R.$ 
 We recall here that $\tau$ is the intensity of the Poisson point process
 ${\mathcal P}.$ 
 Again, 
 using the methods and results of Section 4.3 in \cite{BY} we get
 \begin{proposition}\label{VarProp2}
  For each $f \in {\mathcal C}_b(\R)$ there exists the limit 
  $$ V_f[\xi] := \lim_{\lambda \to \infty} \lambda 
     \Var\left( \langle f, \bar{Z}_{\lambda}^{\xi} \rangle \right) =
     \lim_{\lambda \to \infty} \lambda \Var\left( \langle f, Z_{\lambda}^{\xi} \rangle \right).$$
  Moreover, we have
  \begin{equation}\label{VDEF}
   V_f[\xi] = \tau \langle f \otimes f, \mu \rangle,
  \end{equation}
  where $\mu \in {\mathcal M}(\R \times \R)$ is given by 
  $$
   \mu(A_1 \times A_2) := {\Bbb P}(\xi({\bf 0},{\mathcal P}) \in A_1 \cap A_2)
   + $$
  \begin{equation}\label{MUDEF}
   \frac{1}{2} \tau \int_{\R^d} [{\Bbb P}(\xi({\bf 0},{\mathcal P} \cup x) \in A_1
     ,\; \xi(x,{\mathcal P} \cup {\bf 0}) \in A_2)
    - {\Bbb P}(\xi({\bf 0},{\mathcal P}) \in A_1)
      {\Bbb P}(\xi({\bf 0},{\mathcal P}) \in A_2)] dx
  \end{equation}
  for Borel $A_1, A_2 \subseteq \R^d$ and with $f \otimes f(x,y) :=
  f(x) f(y),\; x,y \in \R.$ The convergence of the integral in
  (\ref{MUDEF}) is guaranteed by the exponential stabilization of $\xi.$ 
 \end{proposition}
 We note that the multiplication rather than division by $\lambda$ in the
 definition of $V_f[\xi]$ above is due to the $\lambda^{-1}$-normalization
 already present in the definition of $Z^{\xi}_{\lambda}.$ 
 The following condition is a natural counterpart of the level-3 condition
 {\bf (L)}. 
 
 {\bf (L')} The log-Laplace functional
   \begin{equation}\label{LaplLev2}
    L^{\xi}_{\lambda;\alpha_{\lambda}}[f] :=  \frac{1}{\alpha_{\lambda}^2} 
    \log {\Bbb E} \exp \left( \langle f, \bar{Z}^{\xi}_{\lambda} \rangle \right)
   \end{equation}
   
   satisfies
   \begin{equation}\label{Var2}
    \lim_{\lambda \to \infty} L^{\xi}_{\lambda;\alpha_{\lambda}}[f] 
    = \frac{1}{2} V_f[\xi].
   \end{equation}

  Under appropriate additional conditions {\bf (L')} would follow as a direct
  consequence of {\bf (L)}, indeed, taking $\hat{\Phi}_f(\xi[\sigma])$ to be
  $\sum_{x \in \sigma \cap [0,1]^d} f(\xi(x,\sigma))$ we see that $\langle 
  f, \bar{Z}^{\xi}_{\lambda} \rangle$ differs from $\langle \hat{\Phi}_f,
  \bar{\Psi}^{\xi}_{\lambda} \rangle$ just by a boundary-order term, which
  can be easily dealt with e.g. by considering a periodised version of the
  process on a torus, thus getting rid of such boundary effects.    
  The point is, though, that thus defined $\hat{\Phi}_f$ is usually not bounded.
  On the other hand, for all our examples both {\bf (L)} and {\bf (L')} do
  follow from the same theory developed in \cite{BESY}. Therefore we have
  decided to formulate both these conditions separately, without resorting
  to tedious general considerations which would not add any extra examples
  to our list of applications.  
 
  The following level-2 moderate deviation theorem is our second main result. 

 \begin{theorem} \label{MDPlevel2}
  For $\xi$ satisfying both the exponential stabilization condition
  {\bf (E)} and the condition {\bf (L')} ,
  the family $\alpha_{\lambda}^{-1} \lambda^{1 \slash 2}
  \bar{Z}_{\lambda}$ satisfies a MDP on ${\mathcal M}(\R)$, 
   endowed with the ${\mathcal C}_b(\R)$-weak topology, with
   speed $\alpha_{\lambda}^2$ and a convex, good rate function
  \begin{equation} \label{rateemp}
   J^{\xi}(\gamma) := \sup_{f \in {\mathcal C}_b(\R)}
    ( \langle f, \gamma \rangle - \frac{\tau}{2} \langle f \otimes f, \mu \rangle ).
  \end{equation}
\end{theorem}

\section{Examples}\label{EXAMPLES}
 Below, we discuss examples of stabilizing functionals for which our 
 general level-3 and level-2 theory applies. Our presentation is
 borrowed from (\cite{BESY}) where level-1 moderate deviation
 principles are established for these functionals. It should be
 noted that the corresponding central limit theorems, under 
 much milder conditions (no homogeneity required) have been
 established in \cite{BY}.

\subsection{Random sequential packing}\label{packing}

 The following
 prototypical random sequential packing model arises in diverse
 disciplines, including physical, chemical, and biological
 processes. See \cite{PY2} for a discussion of the many
 applications, the many references, and also a discussion of
 previous mathematical analysis. In one dimension, this model is
 often referred to as the R\'enyi car parking model \cite{Re}.

 With $N(\tau \la)$ standing for a Poisson random variable with
 parameter $\tau \la$, let $B_{\la,1},$ $B_{\la,2},$ $...,$
 $B_{\la,N(\tau \la)}$ be a sequence of $d$-dimensional balls of
 volume $1$ whose centers are  i.i.d. random $d$-vectors
 $X_1,...,X_{N(\la)}$ uniformly distributed over $Q_{\la} = 
 [-\sqrt[d]{\la} /2, \sqrt[d]{\la} /2]^d.$ 
 Without loss of generality, assume that the balls are sequenced
 in the order determined by marks (time coordinates) in $[0,1]$.
 Let the first ball $B_{\la,1}$ be {\em packed}, and recursively
 for $i=2,3,\ldots, N(\tau\lambda)$, let
 the $i$-th ball $B_{\la,i}$ be  packed iff $B_{\la,i}$ does not
 overlap any ball in $B_{\la,1},...,B_{\la,i-1}$ which has already
 been packed. If not packed, the $i$-th ball is discarded.

 For any finite point set  $\XX \subset \R^d$, assume the points $x
 \in \XX$  have time coordinates which are independent and
 uniformly distributed over the interval $[0,1]$.  Assume unit
 volume balls centered at the points of $\XX$ arrive sequentially
 in an order determined by the time coordinates, and assume as
 before that each ball is packed or discarded according to whether
 or not it overlaps a previously packed ball.  Let $\xi(x; \XX)$ be
 either $1$ or $0$ depending on whether the ball centered at $x$ is
 packed or discarded. Letting $\XX = {\mathcal P}$ we easily see
 that $\xi(\cdot;{\mathcal P})$ describes the random sequential
 packing process as constructed above. This process depends not 
 only on the spatial locations of points but also on their
 $[0,1]$-valued arrival time marks. However, this clearly does
 fit into our general setting by a simple generalisation to the
 marked case.  

 From \cite{BY, PY4} we know that $\xi$ satisfies the exponential
 stabilization condition ${\bf (E)}.$ Moreover, by Section 2 and
 Subsections 6.1 and 6.2 of \cite{BESY} we see that $\xi$ satisfies
 both the {\bf (L)} and {\bf (L')} conditions. In particular,
 our Theorems \ref{GLOWNE} and \ref{MDPlevel2} do apply for
 the random sequential packing functional $\xi.$  To be able to obtain 
{\bf (L)} and {\bf (L')} in \cite{BESY}, we had to apply stabilization methods, 
cumulant techniques, and exponential modification of measures.

\subsection{Spatial birth-growth models}\label{growth}

 Our results for the prototypical packing measures as described
 in Subsection \ref{packing} above, extend to measures arising
 from more general packing models.
 Consider for example the following spatial birth-growth model in
 $\R^d$. Let $\tilde{\mathcal P} := \{(X_i,T_i) \in \R^d \times [0,1] \}$
 be a spatial-temporal Poisson point process.  Seeds appear at uniformly 
 random locations $X_i \in Q_{\la}$ at times $T_i$ i.i.d. and uniform in $[0,1]$.
 When a seed is born, it has initial radius $\rho_i, \ \rho_i \leq L < \infty$,
 and thereafter the radius grows at a constant speed $v_i$, generating a cell
 growing radially in all directions. When one expanding cell touches another,
 they both stop growing in their respective directions. In any event, we assume
 that the seed radii are deterministically bounded, i.e., they never
 exceed a fixed cut-off and they stop growing upon reaching it.
 Moreover, if a seed appears at $X_i$ and if the ball centered at $X_i$ with
 radius $\rho_i$ overlaps any of the existing cells, then the seed
 is discarded. Variants of this well-studied process are used to
 model crystal growth \cite{SKM}.

 To proceed, for any finite point set $\XX \subset \R^d$, assume
 the points $x \in \XX$ have i.i.d. time marks over $[0,1]$.
 A mark at $x \in \X$ represents the arrival time of a seed at $x$.
 Assume that the seeds are centered at the points of $\XX$,
 that they arrive sequentially in an order determined by the
 associated marks, and  that each seed is accepted or rejected
 according to the rules above. Let $\xi(x;\XX)$ be either $1$ or $0$
 according to whether the seed centered at $x$ is accepted or not.
 Letting $\XX = {\mathcal P} \cap Q_{\la}$ we see that
 $\xi(\cdot;{\mathcal P} \cap Q_{\la})$
 corresponds to the spatial birth-growth model introduced above. 

 Again, from \cite{BY, PY4} we know that $\xi$ satisfies the exponential
 stabilization condition ${\bf (E)}.$ Moreover, by Section 3 and
 Subsections 6.1 and 6.2 of \cite{BESY} we see that $\xi$ satisfies
 both the {\bf (L)} and {\bf (L')} conditions. In particular,
 our Theorems \ref{GLOWNE} and \ref{MDPlevel2} do apply for
 the random birth-growth functional $\xi.$

 \begin{remark} 
  The results of the present subsection extend to more general versions
  of the prototypical packing model. The stabilization analysis of
  \cite{PY2} combined with \cite{BESY} yields {\bf (E)} and
  {\bf (L),(L')} in the finite input setting for the number of
  packed balls in the following general models: (a) models with
  balls replaced by particles of random (bounded) size/shape/charge,
 (b) cooperative sequential adsorption models, and (c) ballistic
  deposition models (see \cite{PY2} for a complete description of
 these models). In each case,  our general results apply to the
 functionals $\xi$ putting $1$ in the centers of accepted objects
 and $0$ in the centers of rejected objects.  
\end{remark}

\subsection{Germ-grain models}\label{GGMODEL}
 Let $X_i, \ i \geq 1,$ be i.i.d. uniformly distributed over $Q_{\la}.$
 Let $T, T_i, \ i \geq 1$, be i.i.d. bounded random variables,
 independent of the random variables  $X_i, 1 \leq  i \geq 1.$
 Consider the random grains $X_i + B_{T_i}(\0)$ as well as
 the random set
 $$
  \Xi_{\la} := \bigcup_{i=1}^{N(\tau\la)} (X_i + B_{T_i}(\0)),
 $$
 where $B_r(x)$ again denotes the Euclidean ball centered at $x \in
 {\Bbb R}^d$ of radius $r >0$. The random set $\Xi_{\la}$ usually
 goes under the name of a Boolean model (see e.g. Hall \cite{Ha},
 pp. 141, 233 and Molchanov \cite{Mol} Section 3.2, Example 2.2, p. 35).

 For all $u \in \R^d$, let $T(u)$ be i.i.d. random variables with
 distribution equal to that of $T$. For all $x \in \R^d$ and all
 locally finite point sets $\XX \subset \R^d$, denote by $V(x, \XX)$
 the Voronoi cell around $x$ with respect to $\XX$ and let $\xi(x; \XX)$ be
 the Lebesgue measure of the intersection of $\bigcup_{u \in \XX} B_{T(u)}(u)$
 and $V(x, \X)$.

 For $\xi$ thus defined, we see that $\sum_{x \in {\mathcal P} \cap Q_{\la}}
 \xi(x;{\mathcal P})$ is just the Lebesgue measure of $\Xi_{\la}.$  
 Likewise, we can easily construct a functional $\xi'$ such that
 $\sum_{x \in {\mathcal P} \cap Q_{\la}} \xi'(x;{\mathcal P})$
 coincides with the surface area measure of $\Xi^{\la}$ by defining 
 $\xi'(x;\XX)$ to be the surface area measure of the part of $\partial \Xi_{\la}$
 falling into $V(x,\XX).$ 

 Using \cite{BY} and \cite[Section 6.3]{BESY} we again see that the functionals
 $\xi$ and $\xi'$ as defined above do satisfy both the {\bf (E)}
 and {\bf (L) + (L')} conditions, whence our general results 
 apply. Note that the arguments used when proving conditions
 {\bf (L) + (L')}  for Germ-grain models in \cite{BESY} differ 
 from those used for the packing models, see Section 6.3 ibidem.

\subsection{$k$-nearest neighbors random  graphs}\label{NNGRAPHS}

 Let $k$ be  a positive integer. Given a locally finite  point set
 $\XX \subset \R^d$, the $k$-nearest neighbors (undirected) graph
 on $\XX$, denoted $NG(\XX)$, is the graph with vertex set
 $\XX$ obtained by including $\{x,y\}$ as an edge whenever $y$ is
 one of the $k$ nearest neighbors of $x$ and/or $x$ is one of the
 $k$ nearest neighbors of $y$. The $k$-nearest neighbors (directed)
 graph on $\XX$, denoted $NG'(\XX)$, is the graph with vertex set
 $\XX$ obtained by placing a directed edge between each point and
 its $k$ nearest neighbors.

 For all $t > 0$, let $\xi^t(x; \XX):= 1$ if the length of the edge
 joining $x$ to its nearest neighbor in $\XX$ is less than $t$ and
 zero otherwise. Moreover, for $m \in {\Bbb N}$ we shall consider
 functionals $\xi^{NG}_m$ and $\xi^{NG'}_m$ taking value $1$ if
 the degree of the vertex $x$ in $NG(\XX)$ (respectively $NG'(\XX)$)
 is $m,$ and value $0$ otherwise. Clearly, as usual we shall take
 $\XX := {\mathcal P}.$ It follows now from \cite{BY} and
 \cite[Section 6.3]{BESY} that all the functionals $\xi^t, \xi^{NG}_m$ and
 $\xi^{NG'}_m$ do satisfy {\bf (E)}, {\bf (L)} and ${\bf (L')},$
 whence our general results do apply.

\section{Proof of Theorem \ref{GLOWNE}}
 In view Proposition \ref{VarianceLemma} and condition {\bf (L)} the projective
 limit technique, see Corollary 4.6.11 in \cite{DZ}, allows us to conclude
 that $\bar{\Psi}^{\xi}_{\lambda}$ satisfies the moderate deviation 
 principle in the {\it algebraic dual} $[B_{loc}(\xi[\Sigma])]'$
 endowed with $B_{loc}(\xi[\Sigma])$-weak topology, with the good
 rate function 
 \begin{equation}\label{DZIWNEI}
  [\Lambda^{\xi}]^*(\hat\Theta) := \sup_{\hat{\Phi} \in B_{loc}(\xi[\Sigma])}
   ( \langle \hat{\Phi}, \hat \Theta \rangle - \frac{1}{2} V[\xi;\hat{\Phi}]),\;\;
   \hat\Theta \in [B_{loc}(\xi[\Sigma])]'.
 \end{equation} 
 In view of Theorem \ref{ZASWAR2} below, we have
 \begin{equation}\label{ROWNLA1}
  [\Lambda^{\xi}]^*(\hat\Theta) = I^{\xi}(\hat\Theta|\xi[\Pi])
 \end{equation}
 for $\hat\Theta \in {\mathcal M}^{0,\theta}_{loc}(\xi[\Sigma]).$
 Further, it is easily seen that $\hat{\mathcal T}^{\xi}$ coincides with
 the topology on ${\mathcal M}^{0,\theta}_{loc}(\xi[\Sigma])$ induced
 by the inclusion of this space in $[B_{loc}(\xi[\Sigma])]'$ topologized
 as above. Thus, in view of Lemma 4.1.5 in \cite{DZ}, Theorem \ref{GLOWNE}
 will follow once we show that 
 \begin{equation}\label{ROWNLA2}
  [\Lambda^{\xi}]^*(\hat\Theta) = +\infty
 \end{equation}
 for $\hat\Theta \in [B_{loc}(\xi[\Sigma])]' \setminus {\mathcal M}^{0,\theta}_{loc}(\xi[\Sigma]).$
 To this end, take $\hat\Theta$ with $[\Lambda^{\xi}]^*(\hat\Theta) < +\infty$ and use 
 Proposition \ref{VarianceLemma} writing
 $$ \left\langle \frac{\hat\Phi}{||\hat\Phi ||_{\infty}}, \hat\Theta \right\rangle \leq
    [\Lambda^{\xi}]^*(\hat\Theta) + V\left[\xi;\frac{\hat{\Phi}}{||\hat\Phi ||_{\infty}}\right]
    \leq [\Lambda^{\xi}]^*(\hat\Theta) + V_{D(\hat{\Phi})}[\xi] $$
 for all $\hat\Phi \in B_{loc}(\xi[\Sigma]).$ Consequently, we see that 
 $\hat\Theta$ is a bounded linear form on $\{ \hat{\Phi} \in B_{loc}(\xi[\Sigma]),\;
  \hat \Phi \mbox{ is } {\mathcal F}_{A}-\mbox{measurable } \}$ for each bounded
 Borel $A \subseteq {\Bbb R}^d.$ Using 
 the Riesz representation theorem for the restrictions of $\hat\Theta$ to
 subspaces of functions of $B_{loc}(\xi[\Sigma])$ depending only on the
 marked point configuration within $[-N,N]^d,\; N \to \infty,$ we conclude 
 that $\hat\Theta \in {\mathcal M}_{loc}(\xi[\Sigma]).$ Note that this
 application of the Riesz representation theorem is justified because 
 for each $N \in {\Bbb N}$ the space of finite point configurations 
 in $[-N,N]^d$ can be embedded in the space of compact subsets of $[-N,N]^d$
 endowed with the usual compact Hausdorff metric and with the resulting 
 Borel $\sigma$-field coinciding with ${\mathcal F}_{[-N,N]^d}.$ 
 To complete the proof it is now enough to exclude the case
 $\hat\Theta \in {\mathcal M}_{loc}(\xi[\Sigma]) \setminus {\mathcal M}_{loc}^{0,\theta}(\xi[\Sigma]).$
 However, this is easily done by noting that 
 $$ {\Bbb P}\left(\bar\Psi^{\xi}_{\lambda} \in {\mathcal M}_{loc}(\xi[\Sigma]) \setminus
             {\mathcal M}_{loc}^{0,\theta}(\xi[\Sigma]) \right) = 0 $$
 since $\Psi^{\xi}_{\lambda}$ is translation invariant and has $0$ total mass by
 its definition, and by observing that the space
 ${\mathcal M}_{loc}^{0,\theta}(\xi[\Sigma])$ is closed in ${\mathcal M}_{loc}(\xi[\Sigma])$
 with respect to the $B_{loc}(\xi[\Sigma])$-weak topology.
 The proof is complete. $\Box$

\section{Proof of Theorem \ref{MDPlevel2}}
The proof is organised similar to that of Theorem \ref{GLOWNE}.
In view of Proposition \ref{VarProp2} and the condition {\bf (L')}, 
the projective limit technique, see Corollary 4.6.11 in \cite{DZ},
allows us to conclude that $\bar{Z}^{\xi}_{\lambda}$ satisfies
the moderate deviation principle in the {\it algebraic dual}
$[{\mathcal C}_b(\R)]'$ endowed with ${\mathcal C}_b(\R)$-weak topology,
with the good rate function 
 \begin{equation}\label{POZ2I}
  [L^{\xi}]^*(\gamma) := \sup_{f \in {\mathcal C}_b(\R)} 
   ( \langle f, \gamma \rangle - \frac{1}{2} V_f[\xi]),\;\;
   \gamma \in [{\mathcal C}_b(\R)]'.
 \end{equation} 
In view of Lemma 4.1.5 in \cite{DZ} and of Proposition \ref{VarProp2}
guaranteeing that $[L^{\xi}]^*(\gamma) = J^{\xi}(\gamma)$ for $\gamma
\in {\mathcal M}(\R),$ to complete the proof of Theorem \ref{MDPlevel2}
it is now enough to show that
\begin{equation}\label{InfLev2}
 [L^{\xi}]^*(\gamma) = +\infty,\;\; \gamma \in [{\mathcal C}_b(\R)]' \setminus {\mathcal M}(\R). 
\end{equation} 
To this end, take $\gamma$ with $[L^{\xi}]^*(\gamma) < +\infty$ and
write for $f \in {\mathcal C}_b(\R)$
$$ \left\langle \frac{f}{||f||_{\infty}}, \gamma \right\rangle \leq [L^{\xi}]^*(\gamma) + 
   \frac{1}{2} V_{f/||f||_{\infty}}[\xi] =
   [L^{\xi}]^*(\gamma) + \frac{\tau}{2} \left\langle \frac{f}{||f||_{\infty}} \otimes
   \frac{f}{||f||_{\infty}} , \mu \right\rangle. $$ 
Since the RHS is bounded, this means that $\gamma$ is a bounded operator on ${\mathcal C}_b(\R)$
and hence $\gamma \in {\mathcal M}(\R)$ as required. This completes the proof 
of Theorem \ref{MDPlevel2}. \qed



\section{Properties of the specific relative information}\label{GESTINFWL}
 In this section we discuss a number of properties of the specific
 relative information, as introduced in (\ref{GESTINF}). Our main 
 purpose below is to identify the rate function in Theorem \ref{GLOWNE}. 
 
\subsection{Existence}

\begin{lemma}\label{ISTNIENIE}
 For each translation invariant local null measure $\Theta$
 on $\Sigma$ there exists the limit
 $$ I(\Theta|\Pi) 
    := \frac{1}{2} \lim_{\lambda \to \infty} \frac{1}{\lambda} 
      \int_{\Sigma_{Q_{\lambda}}}
      \left(\frac{d\Theta_{|Q_{\lambda}}}
       {d \Pi_{Q_{\lambda}}}\right)^2
       d \Pi_{Q_{\lambda}}. $$
 Moreover, we have
 \begin{equation}\label{SUPR}
   I(\Theta|\Pi) 
    = \frac{1}{2} \sup_{\lambda \to \infty} \frac{1}{\lambda} 
      \int_{\Sigma_{Q_{\lambda}}}
      \left(\frac{d\Theta_{|Q_{\lambda}}}
       {d \Pi_{Q_{\lambda}}}\right)^2
       d \Pi_{Q_{\lambda}}.
  \end{equation}
 \end{lemma}

\paragraph{Proof}
 For a bounded region $A \subseteq {\Bbb R}$ write
 \begin{equation}\label{IMIEJ}
  I_A(\Theta|\Pi) :=
  \frac{1}{2} \int_{\Sigma_{A}}
      \left(\frac{d\Theta_{|A}}
       {d \Pi_A}\right)^2 d \Pi_A.
 \end{equation}
 It is clear that, by standard superadditivity argument, 
 the proof will be completed once we show that for bounded
 and disjoint $A,B \subseteq {\Bbb R}^d$   
 \begin{equation}\label{NADD}
  I_{A \cup B}(\Theta|\Pi) \geq
  I_{A}(\Theta|\Pi) + I_{B}(\Theta|\Pi). 
 \end{equation}
 To establish (\ref{NADD}), write $\rho_A$ for the density $d\Theta_{|A}
 \slash d \Pi_A,$ define $\rho_B$ and
 $\rho_{A \cup B}$ likewise and let $\rho_{B|A}(\sigma_B|\sigma_A)
 := \rho_{A\cup B}(\sigma_A \cup \sigma_B) - \rho_A(\sigma_A),$
 with $\sigma_A$ and $\sigma_B$ standing
 for generic elements of $\Sigma_A$ and $\Sigma_B$ respectively. It is clear
 that
 \begin{equation}\label{WARUNEK}
  \int_{\Sigma_B} \rho_{A \cup B}(\sigma_A \cup \sigma_B)
  d \Pi_B(\sigma_B) = \rho_A(\sigma_A)
  \mbox{ and hence }
  \int_{\Sigma_B} \rho_{B|A}(\sigma_B|\sigma_A) d \Pi_B(\sigma_B) = 0
  \;\;\; \Pi_A \mbox{ a.s.}
 \end{equation}
 Moreover, since $\int_{\Sigma_A} \rho_A d\Pi_A = \Theta(\Sigma) = 0,$ 
 interchanging $A$ and $B$ in (\ref{WARUNEK}) we are led to
 \begin{equation}\label{WARUNEK2}
  \int_{\Sigma_A} \rho_{A \cup B}(\sigma_A \cup \sigma_B)
  d \Pi_A(\sigma_A) = \rho_B(\sigma_B)
  \mbox{ and }
  \int_{\Sigma_A} \rho_{B|A}(\sigma_B|\sigma_A) d \Pi_A(\sigma_A) = 
  \rho_B(\sigma_B) \;\;\; \Pi_B \mbox{ a.s.}
 \end{equation}
 With this notation we get, using (\ref{WARUNEK}),
 $$ 2 I_{A \cup B}(\Theta|\Pi) =
    \int_{\Sigma_A} \int_{\Sigma_B} \rho^2_{A \cup B}(\sigma_A \cup
    \sigma_B) d\Pi_B(\sigma_B) d\Pi_A(\sigma_A) = $$ $$ 
   \int_{\Sigma_A} \int_{\Sigma_B} (\rho_A(\sigma_A)
   + \rho_{B|A}(\sigma_B|\sigma_A))^2 d\Pi_B(\sigma_B) d\Pi_A(\sigma_A) = $$
 $$ \int_{\Sigma_A} \rho_A^2 d\Pi_A + 2 \int_{\Sigma_A}
    \rho_A(\sigma_A) \int_{\Sigma_B} \rho_{B|A}(\sigma_B|\sigma_A)
    d\Pi_B(\sigma_B)
    d\Pi_A(\sigma_A) + $$ $$ \int_{\Sigma_A} \int_{\Sigma_B}
    \rho^2_{B|A}(\sigma_B|\sigma_A) d\Pi_B(\sigma_B)
    d\Pi_A(\sigma_A) =
 $$
 $$ \int_{\Sigma_A} \rho_A^2 d\Pi_A + \int_{\Sigma_A}
    \int_{\Sigma_B} \rho^2_{B|A}(\sigma_B|\sigma_A)
    d\Pi_B(\sigma_B) d\Pi_A(\sigma_A). $$
 Applying Jensen's inequality we come to
 $$  2 I_{A \cup B}(\Theta|\Pi) \geq  
     \int_{\Sigma_A} \rho_A^2 d\Pi_A + 
     \int_{\Sigma_B}  \left( \int_{\Sigma_A} \rho_{B|A}(\sigma_B|\sigma_A)
     d\Pi_A(\sigma_A) \right)^2 d\Pi_B(\sigma_B). $$
 Using (\ref{WARUNEK2}) we obtain finally
 $$ 2 I_{A \cup B}(\Theta|\Pi) \geq
    \int_{\Sigma_A} \rho_A^2 d\Pi_A + \int_{\Sigma_B} \rho_B^2 d\Pi_B 
   = 2 I_A(\Theta|\Pi) + 2 I_B(\Theta|\Pi) $$
 which yields (\ref{NADD}) and hence completes the proof of 
 the lemma. $\Box$ 

\subsection{Finite volume variational principle and lower semicontinuity}

\begin{lemma}\label{SKZW}
 For a bounded region $A \subseteq {\Bbb R}^d$ we have for each
 $\Phi \in B(\Sigma_A)$
 \begin{equation}\label{JEDNO}
  \frac{1}{2} \Var(\Phi({\mathcal P}_A)) = \sup_{\Theta \in {\mathcal M}^0(\Sigma_A)}
  (\left \langle \Phi, \Theta \rangle - I_A(\Theta|\Pi) \right)
 \end{equation}
 with ${\mathcal P}_A$ standing for the restriction of ${\mathcal P}$ to $A$
 and where ${\mathcal M}^0(\Sigma_A)$ is the collection of all $0$-total mass
 signed measures on $\Sigma_A.$  Moreover, for each
 $\Theta \in {\mathcal M}^0(\Sigma_A)$ we have
 \begin{equation}\label{DRUGIE}
  I_A(\Theta|\Pi) = \sup_{\Phi \in B(\Sigma_A)} \left( \langle 
  \Phi, \Theta \rangle - \frac{1}{2} \Var(\Phi({\mathcal P}_A)) \right).
 \end{equation}
\end{lemma}

\paragraph{Proof}
 Fix $\Phi \in B(\Sigma_A)$ and note that for $\Theta \in {\mathcal M}^0(\Sigma_A)$
 absolutely continuous w.r.t. $\Pi_A$ we have
 $$ \langle \Phi, \Theta \rangle - I_A(\Theta|\Pi) =
    \int_{\Sigma_A} [\Phi - {\Bbb E}\Phi({\mathcal P}_A)] d\Theta 
   - \frac{1}{2} \int_{\Sigma_A} \left( \frac{d \Theta}{d\Pi_A} \right)^2
     d\Pi_A $$
 because $\int_{\Sigma_A} {\Bbb E}\Phi({\mathcal P}_A) d\Theta = 0.$ 
 Consequently, 
 $$ \langle \Phi, \Theta \rangle - I_A(\Theta|\Pi) =
     \int_{\Sigma_A} [\Phi - {\Bbb E}\Phi({\mathcal P}_A)] \frac{d\Theta}{d\Pi_A}
     d\Pi_A  - \frac{1}{2} \int_{\Sigma_A} \left( \frac{d \Theta}{d\Pi_A}
     \right)^2 d\Pi_A \leq $$
 $$ \frac{1}{2} \int_{\Sigma_A}  [\Phi - {\Bbb E}\Phi({\mathcal P}_A)]^2
    d\Pi_A = \frac{1}{2} \Var(\Phi({\mathcal P}_A)), $$
 where the last inequality follows from $f\rho - \frac{1}{2}\rho^2 \leq
 \frac{1}{2} f^2$ for $f := [\Phi - {\Bbb E}\Phi({\mathcal P}_A)]$ and
 $\rho := \frac{d\Theta}{d\Pi_A}.$ Since $\Theta$ was arbitrary
 with $\Theta \ll \Pi_A$ and $I_A(\Theta|\Pi) = +\infty$ for
 $\Theta \not\ll \Pi_A,$ we conclude that
 \begin{equation}\label{N11}
   \Var(\Phi({\mathcal P}_A)) \geq \sup_{\Theta \in {\mathcal M}^0(\Sigma_A)}
   \left( \langle \Phi, \Theta \rangle - I_A(\Theta|\Pi) \right).
 \end{equation}
 To proceed, let $\Theta^{\Phi} \in {\mathcal M}^0(\Sigma_A)$ be
 given by $d\Theta^{\Phi} := [\Phi - {\Bbb E}\Phi({\mathcal P}_A)] d\Pi_A.$
 We have then
 $$ \frac{1}{2} \Var(\Phi({\mathcal P}_A)) = \frac{1}{2} \int_{\Sigma_A} 
    [\Phi - {\Bbb E}\Phi({\mathcal P}_A)]^2 d\Pi_A =
   \int_{\Sigma_A}  [\Phi - {\Bbb E}\Phi({\mathcal P}_A)] d\Theta^{\Phi}
  - \frac{1}{2} \int_{\Sigma_A}
   \left(\frac{d\Theta^{\Phi}}{d\Pi_A}\right)^2 d\Pi_A = $$
 $$ \langle \Phi, \Theta^{\Phi} \rangle - I_A(\Theta^{\Phi}|\Pi). $$
 Combining these equalities with (\ref{N11}) yields now (\ref{JEDNO}).

 The proof of (\ref{DRUGIE}) is analogous. Fix $\Theta \in
 {\mathcal M}^0(\Sigma_A)$ and write for $\Phi \in B(\Sigma_A)$
 $$  \langle \Phi, \Theta \rangle - \Var(\Phi({\mathcal P}_A)) =
    \int_{\Sigma_A} [\Phi - {\Bbb E}\Phi({\mathcal P}_A)] d\Theta 
   - \frac{1}{2} \int_{\Sigma_A}  [\Phi - {\Bbb E}\Phi({\mathcal P}_A)]^2 d\Pi_A = $$
 $$  \int_{\Sigma_A} [\Phi - {\Bbb E}\Phi({\mathcal P}_A)] \frac{d\Theta}{d\Pi_A}
     d\Pi_A  - \frac{1}{2} \int_{\Sigma_A}
     [\Phi - {\Bbb E}\Phi({\mathcal P}_A)]^2 d\Pi_A \leq $$
 $$ \frac{1}{2} \int_{\Sigma_A} 
     \left( \frac{d \Theta}{d\Pi_A} \right)^2 d\Pi_A = \frac{1}{2} 
      I_A(\Theta|\Pi), $$
 where the last inequality follows from $f\rho - \frac{1}{2} f^2 \leq
 \frac{1}{2} \rho^2$ for $f := [\Phi - {\Bbb E}\Phi({\mathcal P}_A)]$ and
 $\rho := \frac{d\Theta}{d\Pi_A}.$ Since $\Phi$ was arbitrary, we see
 that
 \begin{equation}\label{N22}
   I_A(\Theta|\Pi) \geq \sup_{\Phi \in B(\Sigma_A)} \left( \langle 
   \Phi, \Theta \rangle - \frac{1}{2} \Var(\Phi({\mathcal P}_A)) \right).
 \end{equation} 
 To proceed with the proof of the converse inequality observe first
 that if $\Theta \not\ll \Pi_A,$ the expression $\langle \Phi, \Theta
 \rangle - \frac{1}{2} \Var(\Phi({\mathcal P}_A))$ can be made arbitrarily  
 large by adjusting $\Phi$ on a region in $\Sigma_A$ of non-zero
 total variation for $\Theta$ to which $\Pi_A$ assigns zero mass.
 Now, for $\Theta \ll \Pi_A$ let $\Phi^{\Theta} := \frac{d\Theta}{d\Pi_A}.$
 Observe that ${\Bbb E}\Phi^{\Theta}({\mathcal P}_A) = \Theta(\Sigma_A) = 0.$
 Write
 $$ I_A(\Theta|\Pi) = \frac{1}{2} \int_{\Sigma_A} [\Phi^{\Theta}]^2
    d \Pi_A = \int_{\Sigma_A} \Phi^{\Theta} d\Theta -
    \frac{1}{2} {\Bbb E}[\Phi^{\Theta}]^2 = \langle \Phi^{\Theta},\Theta \rangle
    - \frac{1}{2} \Var(\Phi^{\Theta}({\mathcal P}_A)). $$
 Putting this together with (\ref{N22}) yields (\ref{DRUGIE}). 
 This completes the proof of the lemma. $\Box$

\begin{lemma}\label{POLCIAG}
 The mapping
 $$ ({\mathcal M}^{0,\theta}_{loc,\xi}(\Sigma),{\mathcal T}_{\xi})
     \ni \Theta \mapsto I(\Theta|\Pi) $$
 is convex and lower semicontinuous.
\end{lemma}

\paragraph{Proof}
 The convexity follows immediately by the definition of
 $I(\cdot|\Pi)$ in view of the convexity of finite-volume
 functionals $I_{Q_{\lambda}}(\cdot|\Pi).$ Further, the
 variational formula (\ref{DRUGIE}) represents the 
 finite volume functionals $I_{Q_{\lambda}}(\cdot|\Pi)$
 as suprema over $\Phi \in B(Q_{\lambda})$ of 
 ${\mathcal T}_{\xi}$-continuous functionals, which yields
 the ${\mathcal T}_{\xi}$-lower semicontinuity for
 $I_{Q_{\lambda}}(\cdot|\Pi).$ The required 
 ${\mathcal T}_{\xi}$-lower semicontinuity of
 $I(\cdot|\Pi)$ follows now by (\ref{SUPR}).
 $\Box$

\subsection{Infinite volume variational principle}

\begin{theorem}\label{ZASWAR1}
 For each $\hat\Phi \in B_{loc}(\xi[\Sigma])$ we have
 $$ \frac{1}{2} V[\xi;\hat\Phi] = \sup_{\Theta \in {\mathcal M}^{0,\theta}_{loc,\xi}(\Sigma)}
    \left( \langle \hat\Phi, \xi[\Theta] \rangle - I(\Theta|\Pi) \right). $$
\end{theorem} 
 
\paragraph{Proof}
 We claim first that
 \begin{equation}\label{GORA}
  \sup_{\Theta \in {\mathcal M}^{0,\theta}_{loc,\xi}(\Sigma)}
   \left(\langle \hat\Phi, \xi[\Theta] \rangle -
    I(\Theta|\Pi) \right) \leq \frac{1}{2} V[\xi;\hat\Phi].
 \end{equation}
 For each $\Theta \in {\mathcal M}^{0,\theta}_{loc,\xi}(\Sigma)$
 such that $\Theta \ll_{loc} \Pi$ and $I(\Theta|\Pi) < +\infty$ we
 easily conclude from the exponential stabilization assumption
 {\bf (E)} and from the translational invariance of $\Theta$ that 
 $$ \langle \hat\Phi, \xi[\Theta] \rangle - I(\Theta|\Pi) =
    \lim_{\lambda\to\infty}
    \int_{\Sigma_{Q_{\lambda}}} \langle \hat{\Phi},
         \psi^{\xi}_{\lambda}(\sigma_{Q_{\lambda}}) \rangle  
         d\Theta_{|Q_{\lambda}}(\sigma_{Q_{\lambda}}) -
         \frac{1}{2} \lim_{\lambda \to \infty} \frac{1}{\lambda}
    \left(\frac{d\Theta_{|Q_{\lambda}}}{d\Pi_{Q_{\lambda}}} \right)^2
    d\Pi_{Q_{\lambda}}. $$
 Consequently, using that $\Theta$ is a null-measure, we come
 to  
 $$ \langle \hat\Phi, \xi[\Theta] \rangle - I(\Theta|\Pi) =
    \lim_{\lambda\to\infty} \frac{1}{\lambda} 
    \int_{\Sigma_{Q_{\lambda}}} \lambda [\langle \hat{\Phi},
         \psi^{\xi}_{\lambda}(\sigma_{Q_{\lambda}}) \rangle
         - {\Bbb E}\langle \hat{\Phi},\Psi^{\xi}_{\lambda} \rangle]  
    d\Theta_{|Q_{\lambda}}(\sigma_{Q_{\lambda}}) - \frac{1}{2} \lim_{\lambda \to \infty} 
    \frac{1}{\lambda}
    \left(\frac{d\Theta_{|Q_{\lambda}}}{d\Pi_{Q_{\lambda}}} \right)^2
    d\Pi_{Q_{\lambda}} $$ 
 $$ = \lim_{\lambda\to\infty} \frac{1}{\lambda} 
    \int_{\Sigma_{Q_{\lambda}}}
         \lambda [\langle \hat{\Phi},
         \psi^{\xi}_{\lambda}(\sigma_{Q_{\lambda}}) \rangle
         - {\Bbb E}\langle \hat{\Phi},\Psi^{\xi}_{\lambda} \rangle]
         \frac{d\Theta_{|Q_{\lambda}}}{d\Pi_{Q_{\lambda}}}[\sigma_{Q_{\lambda}}] 
         d\Pi_{Q_{\lambda}}(\sigma_{Q_{\lambda}})
   - \frac{1}{2} \lim_{\lambda \to \infty} \frac{1}{\lambda}
    \left(\frac{d\Theta_{|Q_{\lambda}}}{d\Pi_{Q_{\lambda}}} \right)^2
    d\Pi_{Q_{\lambda}} \leq $$
 $$ \frac{1}{2} \lim_{\lambda\to\infty} \frac{1}{\lambda} \int_{Q_{\lambda}} 
    \left(\lambda [ \langle \hat{\Phi}, \psi^{\xi}_{\lambda}(\sigma_{Q_{\lambda}}) \rangle
         - {\Bbb E}\langle \hat{\Phi},\Psi^{\xi}_{\lambda} \rangle] \right)^2 d \Pi_{Q_{\lambda}}
    (\sigma_{Q_{\lambda}}) = \frac{1}{2} \lim_{\lambda\to\infty} \lambda 
    \Var(\langle \hat\Phi,\bar\Psi^{\xi}_{\lambda} \rangle), $$
 where the last inequality comes from $f \rho - \frac{1}{2} \rho^2 \leq
 \frac{1}{2} f^2$ applied for $f := \lambda [ \langle \hat{\Phi},
 \psi^{\xi}_{\lambda} \rangle - {\Bbb E}\langle \hat{\Phi},\Psi^{\xi}_{\lambda}\rangle ]$
 and $\rho := d\Theta_{|Q_{\lambda}} \slash d\Pi_{Q_{\lambda}}.$ 
 Now, in view of Proposition \ref{VarianceLemma} the last limit equals $\frac{1}{2} V[\xi;\hat\Phi].$ 
 Thus, since $\Theta$ was arbitrary with $\Theta \ll_{loc} \Pi$ with $I(\Theta|\Pi) < + \infty$
 and since $I(\Theta'|\Pi) = +\infty$ for $\Theta' \not\ll_{loc} \Pi,$ we conclude the inequality
 (\ref{GORA}) as required. 

 To establish the converse inequality, for each $N \in {\Bbb N}$ construct
 the measure $\Theta^{\Phi \circ \xi}_N$ by partitioning ${\Bbb R}^d$
 into translates $Q_N[i],\; i \in {\Bbb Z}^d$ of the cube $Q_N$
 and setting
 \begin{equation}\label{TH1}
   \Theta^{\Phi \circ \xi}_N := \frac{1}{N} \int_{Q_N}
   \tau_x [\bigoplus_{i \in {\Bbb Z}^d} \Theta^{\Phi \circ \xi}_{N:i}] dx
 \end{equation}
 with
 \begin{equation}\label{TH2} 
  \frac{d\Theta^{\Phi \circ \xi}_{N:i}}{d\Pi_{Q_N[i]}}[\sigma_{Q_N[i]}] :=
     \lambda [\langle \hat{\Phi},\psi^{\xi}_{Q_N[i]}(\sigma_{Q_N[i]})\rangle
   - {\Bbb E}\langle \hat{\Phi},\Psi^{\xi}_N \rangle], \;\;
  \sigma_{Q_N[i]} \subseteq Q_N[i],
 \end{equation}
 where
 $$ \psi^{\xi}_{Q_N[i]}(\sigma_{Q_N[i]}) = \frac{1}{N} \int_{Q_N[i]} \delta_{\xi[\tau_x \Per_N
    (\sigma_{Q_N[i]})]} dx, $$
 see (\ref{OKRESLENIER}), and where $\bigoplus_{i \in {\Bbb Z}^d} \Theta^{\Phi \circ \xi}_{N:i}$
 is given for a cylinder event $S = S_1 \times \ldots \times S_k,\; S_j \in {\mathcal F}_{Q_N[j]}$
 by 
 \begin{equation}\label{OKROPLUS}
  [\bigoplus_{i \in {\Bbb Z}^d} \Theta^{\Phi \circ \xi}_{N:i}](S) =
   \sum_{j=1}^k \Theta^{\phi\circ\xi}_{N:i}(S_j).
 \end{equation}
 Note that this definition is consistent since all $\Theta^{\Phi\circ\xi}_{N:i}$ are
 null-measures (have their total masses $0$). Intuitively speaking, the above construction
 is the counterpart of taking products of probability measures in our null-measure setting.
 Observe also that, by definition, the measure $\Theta^{\Phi\circ\xi}_{N:i}$
 coincides with the translate $\tau_v \Theta^{\Phi\circ\xi}_{N:j}$ where $v$ is the vector
 joining the center of $Q_N[i]$ to the center of $Q_N[j].$ Again, roughly speaking, 
 this construction can be regarded as a null-measure analogue of taking the product law
 of i.i.d. random objects.

 By exponential stabilization {\bf (E)} is clear that
 $\Theta^{\Phi \circ \xi}_N \in {\mathcal M}^{0,\theta}_{loc,\xi}(\Sigma).$ Moreover, by
 the translation invariance of the Poisson point process $\Pi,$ writing $Q^*_N[A] :=
 \bigcup_{Q_N[i] \cap A \neq \emptyset} Q_N[i]$ and $Q^{\partial}_N[A]
 := Q^*_N[A] \setminus A,\; A \subseteq  {\Bbb R}^d,$ in view of
 (\ref{TH1}) and (\ref{TH2}) above we have 
 $$
  [\rho_N]_{\lambda}(\sigma_{Q_{\lambda}}) :=
  \frac{d[\Theta^{\Phi \circ \xi}_{N}]_{¦Q_{\lambda}}}{d\Pi_{Q_{\lambda}}}
  (\sigma_{Q_{\lambda}}) = $$
 \begin{equation}\label{GESTOSC}
  \int_{Q_N} \int_{\Sigma_{Q^{\partial}_N[\tau_x Q_{\lambda}]}}
    \sum_{Q_N[i] \subseteq Q^*_N[\tau_x Q_{\lambda}]}
    \lambda [\langle \hat\Phi, \psi^{\xi}_{Q_N[i]}([\sigma_{Q_{\lambda}} \cup \sigma]_{|Q_N[i]})
             \rangle - {\Bbb E}\langle \hat{\Phi},\Psi^{\xi}_{Q_N} \rangle]
     d\Pi_{Q^{\partial}_N[\tau_x Q_{\lambda}]}(\sigma) dx. 
  \end{equation}
 Using (\ref{GESTOSC}), Proposition \ref{VarianceLemma} and exponential stabilization {\bf (E)},
 as a consequence of the method of \cite{BY} we get
 \begin{equation}\label{WARIANCJA}
  \frac{1}{2} 
  \lim_{N \to \infty} \lim_{\lambda \to \infty} \frac{1}{\lambda} \int_{\Sigma_{Q_{\lambda}}}
  [\rho_N]_{\lambda}^2 d\Pi_{Q_{\lambda}} = \frac{1}{2} V[\xi;\hat\Phi].
 \end{equation}
 Combining (\ref{GESTOSC}) and (\ref{WARIANCJA}) with
 (\ref{TH1}) and (\ref{TH2}) we can write
 $$ \frac{1}{2} V[\xi;\hat\Phi] =  \lim_{N \to \infty} \lim_{\lambda \to \infty} \frac{1}{\lambda}
  \int_{\Sigma_{Q_{\lambda}}}
  [\rho_N]_{\lambda}^2 d\Pi_{Q_{\lambda}} - \frac{1}{2} 
   \lim_{N \to \infty} \lim_{\lambda \to \infty} \frac{1}{\lambda} \int_{\Sigma_{Q_{\lambda}}}
  [\rho_N]_{\lambda}^2 d\Pi_{Q_{\lambda}} = $$
 $$ \lim_{N\to\infty} \lim_{\lambda\to\infty} \frac{1}{\lambda} \int_{\Sigma_{Q_{\lambda}}}
    [\rho_N]_{\lambda} d[\Theta^{\Phi \circ \xi}_N]_{|Q_{\lambda}}
    - \frac{1}{2} \lim_{N\to\infty} \lim_{\lambda \to \infty} \frac{1}{\lambda}
     \int_{\Sigma_{Q_{\lambda}}}
 \left( \frac{d[\Theta^{\Phi \circ \xi}_{N}]_{¦Q_{\lambda}}}
       {d\Pi_{Q_{\lambda}}} \right)^2 d\Pi_{Q_{\lambda}} = $$
 $$ \lim_{N\to\infty} \left( \langle \hat\Phi, \xi[\Theta^{\Phi \circ \xi}_N]
     \rangle - I(\Theta^{\Phi \circ \xi}_N | \Pi) \right). $$
 This implies that
 \begin{equation}\label{DOL}
  \sup_{\Theta \in {\mathcal M}^{0,\theta}_{loc,\xi}(\Sigma)}
   \left( \langle \hat\Phi, \xi[\Theta] \rangle - I(\Theta|\Pi) \rangle \right)
   \geq \frac{1}{2} V[\xi;\hat\Phi].       
 \end{equation}
  Putting (\ref{GORA}) and (\ref{DOL}) together completes the proof  
  of the theorem. $\Box$

\begin{theorem}\label{ZASWAR2}
 For each $\Theta \in {\mathcal M}^{0,\theta}_{loc,\xi}(\Sigma)$ we
 have
 $$ I(\Theta|\Pi) = \sup_{\hat\Phi \in B_{loc}(\xi[\Sigma^{\xi}])}
    \left( \langle \hat\Phi, \xi[\Theta] \rangle - \frac{1}{2} V[\xi;\hat\Phi] \right).
 $$
\end{theorem}

\paragraph{Proof}
 In view of the convexity and lower semicontinuity of $I(\Theta|\Pi)$
 on $({\mathcal M}^{0,\theta}_{loc,\xi}(\Sigma),{\mathcal T}_{\xi}),$
 as stated in Lemma \ref{POLCIAG}, 
 our assertion follows immediately by the Duality Lemma 4.5.8 in \cite{DZ}
 applied for the duality $[\hat\Phi,\Theta ]_{\xi} := \langle \hat\Phi,
 \xi[\Theta] \rangle.$ $\Box$

 \paragraph{Acknowledgements}
  Special thanks are due to Joe E. Yukich, whose ideas and comments have
  motivated us to prepare this paper. Tomasz Schreiber also wishes to
  gratefully acknowledge the support from the Polish Minister of 
  Scientific Research and Information Technology grant 1 P03A 018 28
  (2005-2007).

 \vskip.5cm
Peter Eichelsbacher, Fakult\"at f\"ur Mathematik, Ruhr-Universit\"at Bochum,
NA 3/68, 44780 Bochum, Germany: 
{\em \tt peter.eichelsbacher@ruhr-uni-bochum.de}

 \vskip.5cm

Tomasz Schreiber, Faculty of Mathematics and Computer Science,
Nicholas Copernicus University, Toru\'n, Poland: \ 
{\em \tt tomeks@mat.uni.torun.pl}

\end{document}